\documentclass[12pt]{article}
\usepackage{commands}

\title{Combinatorial Persistent Homology Transform}
\author[1]{Brittany Terese Fasy}
\author[2]{Amit Patel}
\affil[1]{School of Computing and Dept.\ of Mathematical Sciences, Montana State University}
\affil[2]{Department of Mathematics, Colorado State University}
\date{}

\begin{document}

\maketitle

\begin{abstract}
    The combinatorial interpretation of the persistence diagram as a M\"obius inversion was recently shown to be
    functorial. 
    We employ this discovery to
    recast the Persistent Homology Transform of a geometric complex
    as a representation of a cellulation on $\S^n$ to the category of combinatorial persistence diagrams.
    Detailed examples are provided.
    We hope this recasting of the PH transform will allow for the adoption
    of existing methods from algebraic and topological combinatorics to the study of shapes.
    \blfootnote{BTF is supported by the National Science Foundation under grant no.\ DMS~1664858.
    AP is supported by the Leverhulme Trust grant VP2-2021-008.}
\end{abstract}

%%%%%%%%%%%%%%%%%%%%%%%%%%%%%%%%%%%%%%%%%%%%%%
\section{Introduction}\label{sec:intro}
Persistent homology (PH) captures the ``shape'' of a topological space often
arising from data.
It takes as input a filtration of a space, usually paramterized by the reals, and outputs a multiset of points in the
extended plane known as the persistence diagram.
Each point in the persistence diagram represents a
homological ``feature'' of the topological space and its coordinates provide the
interval of parameters for which that feature is present.
Individually, a persistence diagram is a powerful data analysis tool;
however, a single descriptor is often
not rich enough to capture the intricacies of large, complex data.
For example, the family of filtrations giving rise to the same persistence diagram can be arbitrarily large.
One approach for capturing more information about a space is to
consider not just one filtration of the space but a family of filtrations of the space
and therefore a family of persistence diagrams.

One such family of descriptors and the focus of this paper
is the \emph{persistent homology transform}, introduced
initially in~\cite{turner2014persistent}, which takes as input a space
embedded in $\R^{N+1}$ and represents it---or transforms it---into
a family of persistence diagrams parameterized by the set of all directions~$\S^N$.
The idea behind this transform goes back the early 2000s~\cite{agarwal2006extreme} and
was recently generalized to weighted simplicial complexes
in~\cite{wect2020}.
The PH transform is a \emph{complete} representation of the original shape.
That is, no two different shapes have the
same PH transform.
Although $\S^N$ comprises an uncountable set of directions, if $K$ is a nicely
embedded geometric simplicial complex, then it has a finite
representation~\cite{belton2019reconstructing,curry2018directions,fasy2020faithful, micka2020searching}.
The PH transform and other topological transforms have been applied using a sample of
directions~\cite{crawford2019predicting,turner2014persistent, wect2020, betthauser2018topological,
hofer2017constructing, maria2020intrinsic}, yet there is still a substantial gap
between the theory and this practice~\cite{fasy2018challenges}.

\paragraph{Our Contributions}
Consider a space $X$ with an embedding $\phi : X \to \Rspace^{N+1}$.
The PH transform of $\phi$ is the assignment to every
direction $\mu \in \S^N$ the sublevel set persistence diagram of the height function
$\phi_\mu : X \to \Rspace$ in that direction.
As mentioned above, the persistence diagram is the assignment
to every interval of $\Rspace$ a natural number often called its multiplicity or, in our parlance,
a \emph{charge}.
Unfortunately, there is a major downside to this approach.
In practice, $\phi$ is finitely generated; for example, $\phi$ is a PL embedding of a finite simplicial complex.
Its PH transform consists of an uncountable set of
distinct persistence diagrams; however, when ignoring exact birth/death
coordinates and keeping track of only the relations between birth/death coordinates, we observe that there are distinct
places where \emph{combinatorial} changes in the diagrams happen.
We call two persistence diagrams \emph{combinatorially equivalent} when
the diagrams are related by a homeomorphism from~$\R$ to itself.
In this setting, we expect that most directions $\mu \in \Sspace^N$ have an open
neighborhood~$U \subseteq \S^N$ such that the persistence diagrams associated to any two directions~$\nu, \nu' \in U$
are combinatorially~equivalent.

In that light, we take a combinatorial approach to the PH transform by adopting
the combinatorial framework of persistent homology developed by McClearly and Patel \cite{PHFunc}.
In this setting, filtrations are indexed by an abstract totally ordered, finite poset $P$ as opposed
to the real numbers.
The combinatorial persistence diagram is then the assignment to every non-empty (order-theoretic) interval $[a,b] \subseteq P$ a charge.
Given a PL embedding, we show that there is a cellulation of $\S^n$ such that the combinatorial
persistence diagrams associated to any two directions on a common cell are exactly the same.
This cellular
decomposition is the same decomposition observed
in~\cite{belton2019reconstructing,ghrist2018euler,curry2018directions} for the
PH transform and in~\cite{vineyard} for vineyards.
Further, we harness functoriality of McClearly and Patel
to express the PH transform as a representation of the cellulation on $\S^n$ to the category of
combinatorial persistence diagrams and charge-preserving morphisms.
The cellulation is finite and every combinatorial persistence diagram
is finite making our \emph{combinatorial PH transform} a finite object.
We hope that this recasting of the PH transform will allow for the adoption
of existing methods from algebraic and topological combinatorics to the study of shapes.

%%%%%%%%%%%%%%%%%%%%%%%%%%%%%%%%%%%%%%%%%%%%%%

%%%%%%%%%%%%%%%%%%%%%%%%%%%%%%%%%%%%%%%%%%%%%%
\section{Preliminiaries}\label{sec:prelim}
We use this section to establish background and notation.
Section \ref{sec:embed} introduces geometric complexes and their height functions.
Section \ref{sec:pipeline} introduces the combinatorial persistent homology pipeline of McClearly and Patel.
See also \appendref{background} for a discussion of our notation on posets and categories.
As mentioned in Section \ref{sec:intro}, the combinatorial persistence diagram is indexed
by an abstract totally ordered, finite poset as opposed to the traditional persistence diagram that
is indexed by the real numbers.
Since we are departing from traditional persistent homology and focusing
solely on its combinatorial version, we drop the word ``combinatorial''
from our discussion of combinatorial persistent homology.

%%%%%%%%%%%%%%%%%%%%%%%%%%%%%%%%%%%%%%%%%%%%%%
\subsection{Embeddings of Simplicial Complexes}\label{sec:embed}
We start with an introduction to our main object of study: 
geometric complexes and their height functions. 
We adopt notation similar to~\cite[Section 2.1]{Hatcher:478079}.

Fix a simplicial complex $K$, let $K^0 \subseteq K$ be its vertices,
and denote by $|K|$ its underlying space.
For the purposes of this paper, we assume all simplicial complexes are~finite.
Given a function~$\hat \phi : K^0 \to \Rspace^{N+1}$ defined on its vertices, we extend it
to a PL function on the entire underlying space $|K|$ as follows.
For a simplex~$\sigma = [v_0, \cdots, v_k]$, express every point~$x \in | \sigma |$
by its barycentric coordinates~$(x_0, x_1, \cdots, x_k)$.
Then, the linear extension of $\hat \phi$ to $| \sigma |$ is the map
that sends $x$ to the point~$x_0 \hat \phi(v_0) + x_1 \hat \phi(v_1)+ \cdots + x_k \hat \phi(v_k)$.
In this way, the function~$\hat \phi$ extends to a PL function~$\phi : |K| \to \Rspace^{N+1}$
on the entire complex.
If such an extension is injective, we call $\phi$ a \emph{geometric complex}.

We are interested in examining the geometric complex $\phi : |K| \to \Rspace^{N+1}$ from all directions.
Denote by $\Sspace^N := \{ \mu \in \Rspace^{N+1} : ||\mu|| = 1 \}$ the sphere of
all directions in $\R^{N+1}$.
The \emph{height function} on~$\phi$ along a direction $\mu \in \S^N$ is the continuous
function $\phi_\mu : |K| \to \Rspace^{N+1}$ defined as the dot product
$\phi_\mu(x) := \phi(x) \cdot \mu$.  
For $L \subseteq |K|$, we let
$\phi_\mu(L):=\sup_{x\in L} \phi_\mu(x)$.
In particular, if~$L$ is a simplex in
$K$, then $\phi_\mu(L)$ is the maximum height of the vertices that define $L$.

%%%%%%%%%%%%%%%%%%%%%%%%%%%%%%%%%%%%%%%%%%%%%%

%%%%%%%%%%%%%%%%%%%%%%%%%%%%%%%%%%%%%%%%%%%%%%
\subsection{Persistent Homology Pipeline}\label{sec:pipeline}
We now give a brief introduction to the persistent homology pipeline of McCleary and Patel~\cite{PHFunc}
restricted to the special case of combinatorial one-parameter filtrations and to a fixed
simplicial complex $K$.
This pipeline consists of three categories $\Fil(K)$, $\Mon(n)$, and~$\Fnc(n)$ and two functors $\BD_\ast$ and $\M$:
    \[
        \begin{tikzcd}
            \Fil(K) \ar[rrrr, bend left, "\PH_\ast := \M \circ \BD_\ast"]
            \ar[rr, "\BD_\ast"] && \Mon(n) \ar[rr, "\M"] && \Fnc.
        \end{tikzcd}
    \]
We call $\BD_\ast$ the \emph{birth-death functor}, $\M$ the \emph{M\"obius inversion functor},
and the composition~$\PH_\ast$ the \emph{persistent homology functor}.

\paragraph{Filtrations}

Let $\Delta K$ be the poset consisting of subcomplexes of $K$ ordered by inclusion.
Given a finite, totally ordered poset $P$, a
\emph{$P$-filtration} of~$K$ is a functor $\Ffunc : P \to \Delta K$
such that~$\Ffunc(\top) = K$.
Let  $P$ and~$Q$ be finite, totally ordered posets.
A \emph{filtration-preserving morphism} is a triple $(\Ffunc, \Gfunc, \alpha)$
where $\Ffunc : P \to \Delta K$ and $\Gfunc : Q \to \Delta K$ are $P$- and $Q$-filtrations of
$K$, respectively, and $\alpha : P \to Q$
is a bounded monotone function satisfying the following axiom.
For all $a \in Q$, $\Gfunc(a) = \Ffunc(a^\star)$, where~$a^\star := \max \alpha^{-1}[\bot, a]$, i.e., the maximal element in $P$ that maps to $a$:

If $(\Ffunc, \Gfunc, \alpha)$ and $(\Gfunc, \Hfunc, \beta)$ are filtration-preserving
morphisms, then the composition~$(\Ffunc, \Hfunc, \beta \circ \alpha)$
is a filtration-preserving morphism:
    \[
        \begin{tikzcd}
            P \ar[rd, "\Ffunc"] \ar[r, "\alpha"] & Q \ar[d, "\Gfunc"] \ar[r, "\beta"]
            & R \ar[ld, "\Hfunc"] \\
            & \Delta K.
        \end{tikzcd}
    \]

\begin{defn}
    Fix a simplicial complex $K$.
    Let $\Fil(K)$ be the category whose objects are $P$-filtrations of $K$,
    over all finite, totally ordered posets~$P$, and whose morphisms are filtration-preserving morphisms.
    We call $\Fil(K)$ the \define{category of filtrations of $K$}.
\end{defn}

\paragraph{Monotone Integral Functions}
Assume the usual total ordering on the integer $\Z$.
Let~$P$ and $Q$ be two finite, totally ordered posets and let $f : \Int{P} \to
\Z$ and $g :\Int{Q} \to \Z$
be two monotone integral functions on their poset of intervals; see Appendix \ref{sec:background}.
A \emph{monotone-preserving morphism} is a triple $(f, g, \intfcn{\alpha})$, where
$f : \Int{P} \to \Z$ and $g : \Int{Q} \to \Z$ are monotone integral functions and
$\intfcn{\alpha} : \Int{P} \to \Int{Q}$ is a bounded monotone function induced by a bounded monotone
function $\alpha : P \to Q$ satisfying the following axiom.
For all $I \in \Int{Q}$, $g(I) = f(\max  \intfcn{\alpha}^{-1} [\bot, I])$.
In other words, for every interval $I = [a,b]$, $g(I) = f \big( [a^\star,
b^\star] \big)$ where $a^\star = \max f^{-1}[\bot, a]$
and~$b^\star = \max f^{-1}[\bot,b]$.
        \[
	\begin{tikzcd}
            \Int{P} \ar[rd, "f" ] \ar[rr, "\intfcn{\alpha}"] && \Int{Q} \ar[ld, "g"] \\
	& \Z.   &
	\end{tikzcd}
            \]
If $(f, g, \intfcn{\alpha})$ and $(g,h, \intfcn{\beta})$ are two monotone-preserving
morphisms, then the composition~$(f, h, \intfcn{\beta} \circ \intfcn{\alpha})$
is a monotone-preserving morphism:
    \[
        \begin{tikzcd}
            \Int{P} \ar[rd, "f"] \ar[r, "\intfcn{\alpha}"] & \Int{Q} \ar[d, "g"]
            \ar[r, "\intfcn{\beta}"] & \Int{R} \ar[ld, "h"] \\
            & \Z.
        \end{tikzcd}
    \]

\begin{defn}
Let $\Mon$ be the category whose objects are monotone integral functions, over all
finite, totally ordered posets, and whose morphisms are monotone-preserving morphisms.
We call $\Mon$ the \define{category of  monotone functions}.
\end{defn}

\paragraph{Integral Functions}
Let $P$ and $Q$ be finite, totally ordered posets and let $\sigma : \Int{P} \to
\Z$ and~$\tau : \Int{Q} \to \Z$ be two
\emph{integral functions} on their poset of intervals.
Let~$\alpha : P \to Q$ be a bounded poset function such that for
all $I \in \Int{Q}$,
	\begin{equation}
	\label{eq:charge_preserving}
	g(I) = \sum_{J \in  \alpha^{-1}(I)} f(J).
	\end{equation}
Here, if $ \alpha^{-1}(I)$ is empty, then we interpret the sum as $0$.
We call the triple $(\sigma, \tau, \intfcn{\alpha})$ a
\emph{charge-preserving morphism}, where~$\intfcn{\alpha}: \Int{P} \to \Int{Q}$
is the bounded poset function induced by~$\alpha$.
If $(\sigma, \tau, \intfcn{\alpha})$ and $(\tau, \upsilon, \intfcn{\beta})$ are charge-preserving morphisms,
then $(\sigma, \upsilon, \intfcn{\beta} \circ  \intfcn{\alpha})$ is a charge-preserving morphism:
    \[
        \begin{tikzcd}
            \Int{P} \ar[rd, "\sigma"] \ar[r, "\intfcn{\alpha}"] & \Int{Q} \ar[d,
            "\tau"] \ar[r, "\intfcn{\beta}"] & \Int{R} \ar[ld, "\upsilon"] \\
            & \Z.
        \end{tikzcd}
    \]

\begin{defn}
    Let $\Fnc$ be the category whose objects are integral functions over finite,
    totally ordered posets,
    and whose morphisms are charge-preserving morphisms.
    We call $\Fnc$ the \define{category of integral functions}.
\end{defn}

\begin{rmk}
An important observation is that $\Fnc$ does not have arbitrary colimits.
For example, the following two integral functions do not have a colimit.
Let $P=Q$ be totally-ordered posets
on four elements; in particular, we write $P = Q =\{1 < 2 < 3 < 4\}$.
Let~$\sigma \colon \Int{P} \to \Z$ be the function
\[
    \sigma([i,j]) =
    \begin{cases}
        1, \text{ if } i=1,j=3 \\
        1, \text{ if } i=2,j=4 \\
        0, \text{ otherwise},
    \end{cases}
\]
and let
$\tau \colon \Int{Q} \to \Z$ be the~function
\[
    \tau([i,j]) =
    \begin{cases}
        1, \text{ if } i=1,j=4 \\
        1, \text{ if } i=2,j=3 \\
        0, \text{ otherwise}.
    \end{cases}
\]
If a colimit existed, then it would come from
gluing elements in $P$ and gluing elements in $Q$ until there is a common integral function.
This leaves two contenders for the colimit:
\begin{enumerate}
    \item The integral function $\mu : \Int Z \to \Z$ on the poset $Z = \{ a < b < c \}$ that maps $[a,b] \mapsto 1$,
        $[a,c] \mapsto 1$, and everything else to zero. This is obtained by
        mapping $1$ and $2$ in both~$P$ and $Q$ to $a$, $3$ in $P$ and $Q$ to $b$, and $4$ in $P$ and $Q$ to $c$.
    \item The integral function $\nu : \Int Z \to \Z$ on the poset $Z = \{a < b < c\}$ that maps $[a,c] \mapsto 1$,
        $[b,c] \mapsto 1$, and everything else to zero. This is obtained by
        mapping $1$ in both $P$ and~$Q$ to $a$, $2$ in $P$ and $Q$ to $b$, and $3$ and $4$ in $P$ and $Q$ to $c$.
\end{enumerate}
We note that neither of these integral functions can be obtained from the other by gluing elements.
This means that arbitrary colimits do not exist in $\Fnc$.
\end{rmk}

\paragraph{Birth-Death Functor}
Fix a field $\field$.
The $d$-th birth-death functor $\BD_d$ assigns to a filtration $\Ffunc : P \to \Delta K$
the monotone integral function $\BD_d (\Ffunc) : \Int{P} \to \Z$ defined as
follows.
For every interval $[a,b] \in \Int{P}$, where $b \neq \top$, $\BD_d [a,b]$
is the dimension of the $\field$-vector space of~$d$-cycles in the simplicial complex $\Ffunc(a)$
that are~$d$-boundaries in the larger simplicial complex~$\Ffunc(b)$.
For an interval $[a, \top]$, $\BD_d [a, \top]$ is simply the dimension of the
vector space of~$d$-cycles
in $\Ffunc(a)$.
The birth-death functor takes a filtration-preserving morphism~$(\Ffunc, \Gfunc, \alpha)$ to the
monotone-preserving morphism
$\big( \BD_i (\Ffunc),  \BD_d (\Gfunc), \bar \alpha \big)$, which turns out to be a
monotone-preserving morphism.

\paragraph{M\"obius Inversion Functor}
Given a monotone integral function $f : \Int{P} \to \Z$, there is a unique
integral function $\sigma : \Int{P} \to \Z$,
called the M\"obius inversion of $f$, such that
for all~$J \in P$,
        \begin{equation}\label{eq:MI}
            f(J) = \sum_{I \in \Int{P} : I \leq J} \sigma(I).
        \end{equation}
The functor $\M$ assigns to every monotone-integral function its M\"obius inversion.
For every monotone-preserving morphism $(f,g, \intfcn{\alpha})$ in $\Mon$,
$\big( \M\, f, \M\, g, \intfcn{\alpha})$ is a charge-preserving morphism in $\Fnc$.

\begin{defn}
    We call the composition $\PH_\ast :=  \M \circ \BD_\ast :\Fil(K) \to \Fnc$ the \define{persistent homology
    functor}. It assigns to every filtration $P$ in $\Fil(K)$ its persistence
    diagram~$\PH_\ast (F)$
    and to every filtration-preserving morphism $(F, G, \alpha)$ in $\Fil(K)$ the charge-preserving morphism
    $(\PH_\ast(F), \PH_\ast(G), \intfcn{\alpha})$.
\end{defn}

An obvious but important observation is the following proposition, which follows immediately
from the definition of the birth-death function and Equation (\ref{eq:MI}).

\begin{prop}
\label{prop:total_charge}
Let $F$ be a filtration in $\Fil(K)$. Then, the total charge of its $d$-th persistence diagram,
    $\Sigma_{I \in \Int{P}} \PH_d (F)(I)$, is the dimension of the $d$-th cycle space of $K$.
\end{prop}

%%%%%%%%%%%%%%%%%%%%%%%%%%%%%%%%%%%%%%%%%%%%%%

%%%%%%%%%%%%%%%%%%%%%%%%%%%%%%%%%%%%%%%%%%%%%%

%%%%%%%%%%%%%%%%%%%%%%%%%%%%%%%%%%%%%%%%%%%%%%%
%\section{Persistent Homology Pipeline(N) \todo{Delete entire section.}}\label{sec:cosheaf-pipeline}
%\input{body/cosheaf-pipeline}
%%%%%%%%%%%%%%%%%%%%%%%%%%%%%%%%%%%%%%%%%%%%%%

%%%%%%%%%%%%%%%%%%%%%%%%%%%%%%%%%%%%%%%%%%%%%%
\section{Combinatorial Transform for Geometric Complexes}\label{sec:geometric}
Fix a geometric embedding $\phi : |K| \to \Rspace^{N+1}$ of a simplicial complex $K$.

%%%%%%%%%%%%%%%%%%%%%%%%%%%%%%%%%%%%%%%%%%%%%%%%
\subsection{Constructing the Cellulation}
%%%%%%%%%%%%%%%%%%%%%%%%%%%%%%%%%%%%%%%%%%%%%%%%

Let $E := \{ (v, v') : v, v' \in K^0 \}$ be the set of all unordered, distinct pairs
of vertices of~$K$.
For each pair of vertices $e = (v, v') \in E$, there is a unique plane in $\R^{N+1}$ perpendicular
to the vector~$\phi(v)-\phi(v')$ through the origin; we denote that plane $H_e$.
The intersection of~$H_e$ with the $N$-sphere,~$S_e := H_e \cap \Sspace^{N}$,
is a great~$(N-1)$-sphere of
$\S^n$.
By the Jordan-Brouwer separation theorem,~$\Sspace^{N} \setminus S_e$
has two connected components (often called \emph{sides} or \emph{hemispheres})
denoted arbitrarily by $S_e^+$ and~$S_e^-$;
see~\cite{lebesgue1911,brouwer1911,Hatcher:478079}.
Let~$\Acal := \{ S_e : e \in E\}$ be the multiset of all such great
$N$-spheres.
We call $\Acal$ along with an assignment of a sign to each side of $S_e$, over all $e \in E$,
a \emph{signed arrangement} of spheres.
The signed arrangement $\Acal$ is \emph{essential}
if the following condition holds:~$\bigcap_{e \in \Acal} S_e = \emptyset$.

The signed arrangement $\Acal$ induces a cell complex on $\S^{N}$ where
every cell is described by whether it lies in $S_e^-$, $S_e$, or $S_e^+$
over every $S_e \in \Acal$.
To define this cell complex, let~$\{ -, 0, +\}^E$ denote the set
of all functions from $E$ to $\{-, 0, + \}$.
Consider the following partial order on~$\{ -, 0, + \}$: $0 < -$ and $0 < +$.
This partial order extends to a partial order on the set of functions
$\{ -, 0,+ \}^E$, where~$f \leq g$ if~$f(e) \leq g(e)$, for all $e \in E$.
Let $\Phi : \Sspace^{N} \to \{ -, 0, + \}^E$ be the function that assigns to every direction
$\mu \in \Sspace^N$ the following vector indexed by elements~$e \in E$:
        \begin{equation}\label{eq:vector}
	\Phi(\mu)_e := \begin{cases}
	- & \text{if $\mu \in S_e^-$} \\
	0 & \text{if $\mu \in S_e$} \\
	+ & \text{if $\mu \in S_e^+$}.
	\end{cases}
	\end{equation}
Consider the subposet $\Pcal(\Acal) := \{ \Phi(\mu) : \mu \in \S^N \}
\subseteq \{ - , 0, + \}^E$.
The following proposition is a rewording of \cite[Proposition 5.1.5]{oriented_matroids}.

\begin{prop}
\label{prop:cellulation}
    Let $\phi : |K| \to \Rspace^{N+1}$ be geometric embedding and let
    $\Acal$ be an induced signed arrangement of spheres on $\S^N$.
    If $\Acal$ is essential, then the sets~$\Phi^{-1}(f) \subseteq \S^N$, over
    all~$f \in \Pcal(\Acal)$, are the cells of a cell complex
    $(\S^N , \Sstrat)$.
    Furthermore, for every pair of such cells~$C_1 := \Phi^{-1}(f)$ and $C_2 := \Phi^{-1}(g)$ in $\Sstrat$,
    $C_1 \leq C_2$ if and only if $f \leq g$ in~$\Pcal(\Acal)$.
\end{prop}

The signed arrangement $\Acal$ associated to the piecewise linear embedding $\phi$
is essential if and only if $\bigcap_{e \in \Acal} S_e \neq \emptyset$ if and only if
there is not a direction $\mu \in \S^n$ such that the height
function~$\phi_\mu$ is constant on all of $K$.
For example, if $\phi$ embedds $K$ into a linear $N$-subspace of~$\Rspace^{N+1}$,
then $\Acal$ is not essential and we cannot apply Proposition \ref{prop:cellulation}.
However, we can easily fix this problem by adding one more $(N-1)$-great sphere to
the arrangement that is different from any of the $(N-1)$-great spheres already in $\Acal$.
In general, for $\dim\left( \bigcap_{e \in \Acal} S_e \right) = k$, we have to
add $k+1$ great spheres to the
arrangement $\Acal$.

%%%%%%%%%%%%%%%%%%%%%%%%%%%%%%%%%%%%%%%%%%%%%%%%
\subsection{Filtrations over Cellulation}
%%%%%%%%%%%%%%%%%%%%%%%%%%%%%%%%%%%%%%%%%%%%%%%%
Fix a geometric embedding $\phi : |K| \to \Rspace^N$ and
let~$(\S^N, \Sstrat)$ be the induced cellulation on the~$N$-sphere as constructed above.
Recall $\Sstrat$ is a poset where $C_1 \leq C_2$ if $C_1$ is a face of~$C_2$.
For every direction $\mu \in \S^n$, the height function $\phi_\mu$ gives rise to an object of $\Fil(K)$.
By construction of the cellulation, any two directions $\mu_1$ and $\mu_2$ give rise to the
same object in $\Fil(K)$.
Further, for every face relation $C_1 \leq C_2$, there is a natural filtration-preserving
morphism from the filtration associated to $C_2$ to the filtration associated to $C_1$.
We now formalize this data as a functor $\Ffunc : \Sstrat \to \Fil(K)$.

We start by defining $\Ffunc$ on the cells of $\Sstrat$.
Choose a cell $C \in \Sstrat$ and a direction $\mu \in C$.
Two vertices $v, v' \in K^0$ are related, denoted $v \sim_\mu v'$, if
$\phi_\mu(v) = \phi_\mu(v')$.
The relation~$\sim_\mu$ is an equivalence relation.
Note that $v \sim_\mu v'$ if and only if $\mu \in S_{(v, v')}$.
Denote by~$\Pcat_\mu$ the set of equivalence classes~$K^0 / \sim_\mu$ union the singleton $\{ \top \}$.
For two equivalence classes~$[v]_\mu, [v']_\mu$ in~$\Pcat_\mu$, let
$[v]_\mu \leq [v']_\mu$ if~$\phi_\mu(v) \leq \phi_\mu(v')$.
Make $\top$ the top element.
Thus~$\Pcat_\mu$ is a finite, totally ordered lattice that is independent of the choice of
$\mu \in C$.
Note that $[v]_\mu \neq [v']_\mu$ if and only if $\mu \in S_{(v, v')}^+$
or $\mu \in S_{(v, v')}^-$.
We now define the filtration~$\Ffunc(C) := \Ffunc_\mu : P_\mu \to \Fil$.
For every~$[v]_\mu$, let
\begin{equation}\label{eq:dir-filt}
\Ffunc_\mu\big( [v]_\mu \big) := \big \{ \sigma \in K : \phi_\mu ( \sigma ) \leq \phi_\mu(v) \big \}.
\end{equation}
Let $\Ffunc_\mu(\top)$ be the entire simplicial complex $K$.
Note that $\Ffunc_\mu \big( [v]_\mu \big)$ is a subcomplex of $K$ and for $[v]_\mu \leq [v']_\mu$,
$\Ffunc_\mu \big( [v]_\mu \big)$ is a subcomplex of~$\Ffunc_\mu \big( [v']_\mu \big)$.
In other words, $\Ffunc_\mu$ is an object of~$\Fil(K)$.

We now define $\Ffunc$ on the face relations $C_1 \leq C_2$ of $\Sstrat$.
Choose directions $v \in C_1$ and~$\mu \in C_2$.
There is a canonical bounded lattice function
$\alpha : \Pcat_\mu \to \Pcat_\nu$ that turns out to be a filtration-preserving morphism
$\big(\Ffunc(C_2), \Ffunc(C_1), \alpha \big)$ as follows.
Since $\Phi(v) \leq \Phi(\mu)$,~$[v]_\mu \subseteq [v]_\nu$ for every vertex
$v \in K^0$.
Let $\alpha \big( [v]_\mu \big) := [v]_\nu$.
Again, since $\Phi_{(v, v')}(\mu) \leq \Phi_{(v, v')}(v)$, we have
that~$[v]_\mu \leq [v']_\mu$ implies $[v]_\nu \leq [v']_\nu$.
Let $\alpha (\top) := \top$.
Thus $\alpha$ is a bounded lattice function.
By piecewise linearity of $\phi$,
a simplex~$\sigma \in \Ffunc_\mu \big( [v]_\mu \big)$ if and only if
$\phi_\mu (v_i) \leq \phi_\mu (v)$ for every vertex~$v_i$ of $\sigma$.
Since~$\alpha$ is order-preserving, $\phi_\mu (v_i) \leq \phi_\mu (v)$ implies
$\phi_\nu (v_i) \leq \phi_\nu(v)$ and
thus~$\sigma \in \Ffunc_\nu \big( [v]_\nu \big)$.
Thus~$\big( \Ffunc_\mu, \Ffunc_\nu, \alpha \big)$ is a filtration-preserving morphism.
Now consider a third cell $C_3$ such that~$C_1 \leq C_2 \leq C_3$ and choose a direction
$w \in C_3$.
Let $\beta : \Pcat_w \to \Pcat_\mu$ be the corresponding bounded lattice function.
Then, the assignment $[v]_w \mapsto [v]_\mu$ is the composition
$\alpha \circ \beta \big( [v]_w \big)$.
Thus the composition $\Ffunc(C_1 \leq C_2) \circ \Ffunc(C_2 \leq C_3)$ of filtration-preserving morphisms
equals~$\Ffunc(C_1 \leq C_3)$.

%%%%%%%%%%%%%%%%%%%%%%%%%%%%%%%%%%%%%%%%%%%%%%%%
\subsection{Persistence Diagrams over Cellulation}
%%%%%%%%%%%%%%%%%%%%%%%%%%%%%%%%%%%%%%%%%%%%%%%%
Finally, we define the combinatorial PH transform associated to a geometric complex.

\begin{defn}
\label{defn:transform_cosheaf}
    Let $\phi : |K| \to \Rspace^{N+1}$ be a geometric complex,
    $(\Sspace^N, \Sstrat)$ the induced cellulation on the $N$-sphere, and
    $\Ffunc : \Sstrat \to \Fil$ the functor as constructed above.
    The \define{$d$-th combinatorial PH transform} of $\phi$ is the functor~$\PH_d \circ \Ffunc : \Sstrat \to \Fnc$.
    The \define{combinatorial PH transform} of $\phi$ is the indexed set $\big\{ \PH_d \circ F \big\}_{d \in \N}$,
    over all dimensions~$d$.
\end{defn}

Fix an embedding $\phi : |K| \to \Rspace^{N+1}$ and let $\tilde F_d : \Sstrat \to \Fnc$
be its $d$-th combinatorial PH transform.
The \emph{display locale} \cite{curry, CRG} of $\tilde F_d$, or generalized vineyard \cite{vineyard},
has an interesting structure.
By Proposition \ref{prop:total_charge}, there is a constant $n \in \N$ such
that, for each cell $C \in \Sstrat$, the total charge of the
persistence diagram~$\tilde F_d (C)$ is $n$.
Further, for every cell relation~$C_1 \leq C_2$, the morphism $\tilde F(C_1 \leq C_2)$
taking the persistence diagram $\tilde F(C_2)$ to the persistence diagram~$\tilde F(C_1)$
is charge-preserving.
Thus, there is a unique charge associated to each connected component of the display locale.

\later{
%%%%%%%%%%%%%%%%%%%%%%%%%%%%%%%%%%%%%%%%%%%%%%%%
\subsection{To Remove}
%%%%%%%%%%%%%%%%%%%%%%%%%%%%%%%%%%%%%%%%%%%%%%%%

For a fixed cell complex $(\S^N, \Sstrat)$, denote by $[ \Sstrat, \Fnc]$
the category of functors $\Acat: \Sstrat \to \Fnc$ and natural transformations.
Given two cell complexes $(\S^N, \Sstrat)$ and $(\S^N, \mathcal{D})$ and two
cellular cosheaves $\Acat : \Sstrat \to \Fnc$ and $\Bfunc : \mathcal{D} \to \Fnc$,
we say $\Afunc$ is \emph{isomorphic} to $\Bfunc$ if the following holds:
there is a homeomorphism $h : \S^N \to \S^N$ inducing a poset isomorphism
$h_\ast : \Sstrat \to \mathcal{D}$ such that $\Bfunc \circ h_\ast$ is isomorphic $\Afunc$ in $[ \mathcal{C}, \Fnc]$.
\brittany{maybe add equivalence for our PH transform cosheaves too?  Then,
below, we don't need $d$.}
We say that two PH transform cosheaves $A= \{ A_d\}_{d\in \Z}$ and $B= \{B_d
\}_{d\in \Z}$ are isomorphic if there exists a homeomorphism $h : \S^N \to \S^N$ inducing a poset isomorphism
$h_\ast : \Sstrat \to \mathcal{D}$ such that $B_d \circ h_\ast$ is
isomorphic $A_d$ in $[ \mathcal{C}, \Fnc]$ for all $d \in \Z$.

\begin{thm}[Uniqueness]\label{thm:uniqueness}
    Let  $\phi : |K| \to \Rspace^{N+1}$ and $\psi : |L| \to \Rspace^{N+1}$ be
    piecewise linear embeddings of finite simplicial complexes that take their
    vertex sets to general position.
    Let~$(\S^N, \Sstrat)$ be the cellulation induced by $\phi$, $(\S^N, \mathcal{D})$
    the cellulation induced by $\psi$, $\Afunc_d : \Sstrat \to \Fnc$
    the $d$-th PH transform cosheaf of $\phi$, and $\Bfunc_d : \mathcal{D} \to \Fnc$
    the $d$-th PH transform cosheaf of $\psi$.
    If $\Afunc_d$ is isomorphic to $\Bfunc_d$, for all $d \in \Nspace$, then
    $K$ and $L$ are isomorphic as simplicial complexes.
\end{thm}

\todo{We need to talk through this. Your
paper does not use the language of cosheaves. How can you have a proof?}
The proof follows from the following observations:
that the $N$-cells correspond to the sets of directions where all vertices are
seen at distinct heights (\rmkref{highd-cells}),
that the $i$-simplices are in one-to-one correspondence with birth coordinates
in $\M \circ \BD_i \circ \Ffunc$ and death coordinates in $\M \circ \BD_{i-1} \circ
\Ffunc$,
and
that we can test for the existence of any simplex
by knowing only the partial order of the
vertices in a finite set of given directions (see \cite[Remark 21]{fasy2020faithful}).

}

%%%%%%%%%%%%%%%%%%%%%%%%%%%%%%%%%%%%%%%%%%%%%%

%%%%%%%%%%%%%%%%%%%%%%%%%%%%%%%%%%%%%%%%%%%%%%
%\section{Main Theorem}\label{sec:main-theorem}
%\input{body/main-theorem}
%%%%%%%%%%%%%%%%%%%%%%%%%%%%%%%%%%%%%%%%%%%%%%

\newpage

%%%%%%%%%%%%%%%%%%%%%%%%%%%%%%%%%%%%%%%%%%%%%%
\section{Examples}\label{sec:examples}
In this section, we walk through two examples, one in $\R^2$ and one in $\R^3$.

\subsection{The V Example}

\begin{figure}
    \centering
    \centering
    \includegraphics[height=.75in]{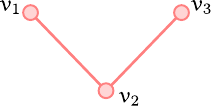}
    \caption{The V, embedded in $\R^2$.  This simplicial complex has three
    vertices and two edges. By exploring the combinatorial PH transform for this
    example, we illustrate each step of the construction.}
    \label{fig:ex-v-gsc}
\end{figure}

We start with a simple example to step through the construction of the combinatorial PH
transform of
\defnref{transform_cosheaf}.
Specifically, our abstract simplicial complex consists of two edges connected along a single
vertex:
$K= \{ \{v_1\}, \{v_2\}, \{v_3\} , \{v_1,v_2\}, \{ v_2,v_3\} \}$.
This
complex is piecewise linearly embedded into $\R^2$ through the function $\phi
\colon |K| \to \R^{2}$ with~$\phi(v_1)=(-1,1)$, $\phi(v_2)= (0,0)$,
and~$\phi(v_3)=(1,1)$. See \figref{ex-v-gsc}.

\begin{figure}
    \centering
    \begin{subfigure}[b]{0.4\textwidth}
        \centering
        \includegraphics[height=2in]{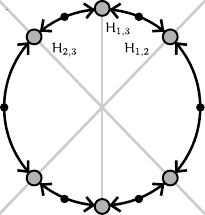}
        \caption{Lines Intersecting $\S^1$}
        \label{fig:ex-v-lines}
    \end{subfigure}
    \hspace{3em}
    \begin{subfigure}[b]{0.4\textwidth}
        \centering
        \includegraphics[height=2in]{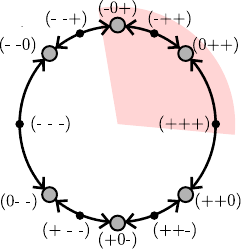}
        \caption{Cellulation of $\S^1$}
        \label{fig:ex-v-sphere}
    \end{subfigure}
    \caption{
    In (\subref{fig:ex-v-lines}), we see the three linear subpaces of $\R^2$ that
    are used to define the cellulation over~$\S^1$.
    In (\subref{fig:ex-v-sphere}), each cell is labeled by a vector in
    $\{-,0,+\}^3$ according to which side of $S_{1,2}$, $S_{1,3}$, and $S_{2,3}$
    the cell falls.
    For example, the vector $(-++)$ labels the one-cell whose points are all
    in~$S_{1,2}^-$,~$S_{1,3}^+$, and~$S_{2,3}^+$.
    The vector $(0++)$ labels the
    zero-cell that is
    in~$S_{1,2}$,~$S_{1,3}^+$, and~$S_{2,3}^+$. In fact,~$\cell_{(0++)}=S_{1,2}\cap S_{1,3}^+\cap S_{2,3}^+$.
    Note that no label is $(000)$, and that all labels are distinct.
    The partial order of the cells is denoted by arrows (where $a \to b$
    indicates that $b < a$).
    }\label{fig:ex-v}
\end{figure}

\paragraph{Cellulation}
The embedding $\phi$ induces a cellulation on the unit circle $\S^1$ of
directions in~$\R^2$.
We now walk through the construction of this cellulation.
Let~$E$ be the set of unique pairs of vertices in $K$.
In other words, let
$E := \{ (v_1,v_2), (v_1,v_3),  (v_2,v_3) \}$.
For every~$(v_i, v_j) \in E$, let~$H_{(i,j)}$ be the unique line through the
origin that is perpendicular to the vector~$\phi(v_i)-\phi(v_j)$.  That is,
$H_{(i,j)}$ is the line
perpendicular to the line through the embeddings of $v_i$ and~$v_j$, as shown in
\subfigref{ex-v}{lines}.
The intersection $S_{(i,j)} := H_{(i,j)} \cap \S^1$ is a zero-sphere separating $\S^1$ into
two components.
As an arbitrary choice of signs,
we denote the connected component containing $(1,0)$ as $S_{(i,j)}^+$ and the other
connected component as~$S_{(i,j)}^-$.
We now have a function~$\Phi \colon \S^1 \to \{-,0,+\}^E$ as in
\eqnref{vector}.
Finally, let~$(\S^1,\Sstrat)$ be the cell complex as defined in
\propref{cellulation}; in particular, $\Sstrat$ has six
zero-cells and six one-cells.
See \subfigref{ex-v}{sphere}, where the $12$ cells
are labeled with their images under~$\Phi$. Notice that in this example, we have
an essential arrangement.
Since no two cells have the same label,
we use $\cell_X$ to denote the cell of $\Sstrat$ labeled $X \in \{-,0,+\}^3$ (if
such a cell exists).
Again, by \propref{cellulation}, we also have a partial order on the
cells~(e.g.,~$\cell_{(0++)} < \cell_{(+++)}$).

\paragraph{Objects of Filtration Functor}
The geometric embedding $\phi$ induces a functor~$\Ffunc : \Sstrat \to \Fil$ of filtrations
and filtration-preserving morphisms.
We now walk through the construction of~$\Ffunc$ on the four cells
in the highlighted region of \subfigref{ex-v}{sphere}.
\begin{figure}
    \centering
    \centering
    \includegraphics[height=1.5in]{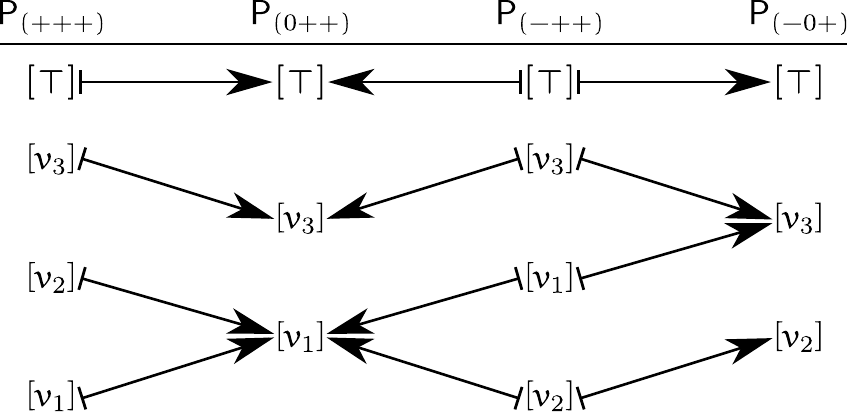}
    \caption{The bounded lattice functions for the three highlighted face relations in
    \subfigref{ex-v}{sphere}.  Notice that the two maps into $\Pcat_{(0++)}$
    are nearly bijections, except for two vertices ($v_1$ and~$v_2$) that map to the
    same equivalence class for both maps. This corresponds to the transposition of the
    two vertices
    between directions in~$\cell_{(+++)}$ and $\cell_{(-++)}$.
    In fact, this property holds
    more generally for any face relation between codimension-one cells.
    } \label{fig:ex-v-lattice}
\end{figure}

Let $\mu$ be a direction in $C_{(+++)}$.
Using the height function $\phi_\mu \colon |K| \to \R$,
we obtain the following totally ordered
lattice $[v_1]_\mu<[v_2]_\mu<[v_3]_\mu$; see \figref{ex-v-lattice} (far left). To this, we add an additional equivalence
class $[\top]_\mu$ as the top element of the lattice; in other words, we have the
totally ordered lattice: \[ [v_1]_\mu<[v_2]_\mu<[v_3]_\mu < [\top]_\mu.\] Denote this lattice~$\Pcat_\mu$.
Associated to $\Pcat_\mu$ is the filtration
$\Ffunc_\mu$ of \eqnref{dir-filt}.  Specifically,
	\begin{align*}
	\Ffunc_\mu([v_1]_\mu) &:= \{v_1\} \\
	\Ffunc_\mu([v_2]_\mu) &:= \{ v_1, v_2, \{v_1,v_2\} \} \\
	\Ffunc_\mu([v_3]_\mu) &:= \{ v_1, v_2, v_3 , \{v_1,v_2\}, \{v_2,v_3\} \} = K \\
	\Ffunc_\mu([\top]_\mu) &:=K.
	\end{align*}
Notice that $\Ffunc_\mu([v_1]_\mu) \subset \Ffunc_\mu([v_1]_\mu) \subset \Ffunc_\mu([v_1]_\mu) \subseteq \Ffunc_\mu([\top]_\mu)$.
Moreover, this lattice and filtration are independent of the
choice of $\mu$.  Therefore, we can define~$[v_i]_{(+++)}:=[v_i]_\mu$,
$[\top]_{(+++)}:=[\top]_\mu$,
$\Pcat_{(+++)}:=\Pcat_\mu$, and $\Ffunc_{(+++)}:=\Ffunc_\mu$.

Next, consider the cell~$\cell_{(-++)}$.  We proceed as above and observe that
$\Pcat_{(-++)}$ is the totally ordered lattice~$[v_2]_{(-++)} < [v_1]_{(-++)} <
[v_3]_{(-++)} < [\top]_{(-++)}$; see \figref{ex-v-lattice} (middle right).
For simplicity, we drop the subscripts of the equivalence
classes when the cell (or direction) is clear from context and simply write
$[v_2]$ from here on.
Again, we use~\eqnref{dir-filt} to define the
filtration $\Ffunc_{(-++)}$ is defined by
	\begin{align*}
	\Ffunc_{(-++)}([v_2]) &:= \{ v_2 \} \\
	\Ffunc_{(-++)}([v_1]) &:=\{ v_1, v_2, \{v_1,v_2\} \} \\
	\Ffunc_{(-++)}([v_3]) &:= \{ v_1, v_2, v_3 , \{v_1,v_2\}, \{v_2,v_3\} \} = K \\
	\Ffunc_{(-++)}([\top]) &:= K.
	\end{align*}
We note the partial order on $\Pcat_{(-++)}$ given by the heights of vertices in
direction $\mu$
induces the following partial order
on the simplicial complexes:~$\Ffunc_{(-++)}([v_2]) \subset \Ffunc_{(-++)}([v_1]) \subset \Ffunc_{(-++)}([v_3]) \subseteq
\Ffunc_{(-++)}([\top])$.

The zero-cells of $\Sstrat$ are exactly the directions that two vertices
are seen at the same height.  For example, in $\cell_{(0++)}$,the heights of
$v_1$ and $v_2$ are the same (since~$\phi_\mu(v_1)=\phi_\mu(v_2)$), and this height is less than $\phi_\mu(v_3)$.  Thus, the
lattice $\Pcat_{(0++)}$ is the following total order on the induced equivalence
classes:  $[v_1] < [v_3] < [\top]$.  In addition, the filtration $\Ffunc_{(0++)}$ is
defined by
	\begin{align*}
	\Ffunc_{(0++)} ([v_1]) &:= \{ v_1, v_2, \{v_1,v_2\} \} \\
	\Ffunc_{(0++)}([v_3]) &:= \{ v_1, v_2, v_3 , \{v_1,v_2\}, \{v_2,v_3\} \} = K \\
	\Ffunc_{(0++)}([\top]) &:= K.
	\end{align*}
Again, we notice the inclusion of subcomplexes: $\Ffunc_{(0++)}([v_1]) \subset
\Ffunc_{(0++)}([v_3]) \subseteq \Ffunc_{(0++)}([\top])$.
Finally, $\Pcat_{(-0+)}$ is the lattice~$[v_2] < [v_3] < [\top]$, with
corresponding
filtration $\Ffunc_{(-0+)}$ defined by~$\Ffunc_{(-0+)} ([v_2]) = \{ v_2\}$ and
$\Ffunc_{(-0+)}([v_3])= \{ v_1, v_2, v_3 , \{v_1,v_2\}, \{v_2,v_3\} \} = K = \Ffunc_{(-0+)}([\top])$.  See \figref{ex-v-lattice}.

\paragraph{Arrows of the Filtration Functor}
We now walk through the construction of $\Ffunc$ on the three arrows
in the highlighted region of \subfigref{ex-v}{sphere}.
Consider the face relation
$\cell_{(0++)} \leq \cell_{(+++)}$.
We define the filtration-preserving morphism
$\big( \Ffunc(\cell_{(+++)}), \Ffunc(\cell_{(0++)}, \alpha \big)$ as follows.
For each vertex $v_i$, we have the
inclusion~$[v_i]_{(+++)} \subseteq [v_i]_{(0++)}$. Indeed,
\begin{align*}
    [v_1]_{(+++)} &= \{v_1\} \subseteq \{v_1,v_2\} = [v_1]_{(0++)}\\
    [v_2]_{(+++)} &= \{v_2\} \subseteq \{v_1,v_2\} = [v_2]_{(0++)}\\
    [v_3]_{(+++)} &= \{v_3\} \subseteq \{v_3\} = [v_3]_{(0++)}\\
    [\top]_{(+++)} &= \emptyset \subseteq \emptyset = [\top]_{(0++)}\\
\end{align*}
As a consequence, $F([v]_{(+++)}))\subset \Ffunc([v]_{(0++)})$ for all objects
$[v] \in
\Pcat_{(+++)}$.
Thus, we can define~$\alpha(\Ffunc([v]_{(+++)})):=\Ffunc([v]_{(0++)}) $.  In other words, we have
collapsed two equivalence classes in the lattice and consequently the filtration,
and left the rest untouched. See \figref{ex-v-lattice}.

\begin{figure}
    \centering
    \begin{subfigure}[b]{0.4\textwidth}
        \centering
        \includegraphics[height=2.5in]{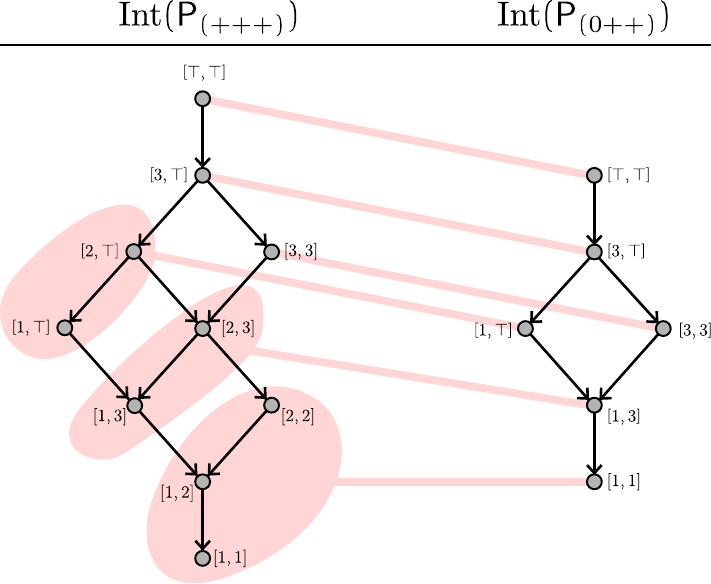}
        \caption{Intervals}
        \label{fig:ex-v-int2zb-interval}
    \end{subfigure}
    \hspace{3em}
    \begin{subfigure}[b]{0.4\textwidth}
        \centering
        \includegraphics[height=2.5in]{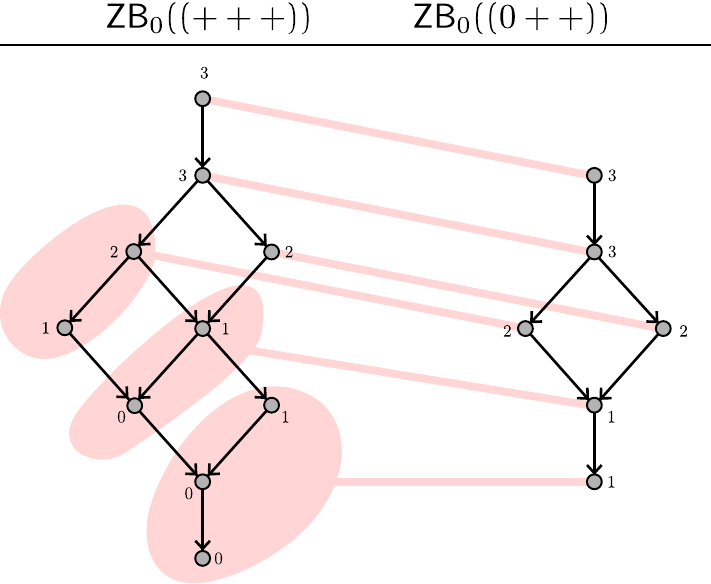}
        \caption{Image of  $\BD_0$}
        \label{fig:ex-v-int2zb-zb}
    \end{subfigure}
    \caption{The maps $\Int{P_{(+++)}} \to \Int{P_{0++}}$ and $\BD_0(\Int{P_{(+++)}} \to
    \Int{P_{0++}})$. In both maps, the objects in the same pink region get mapped to the same
    object in the codomain.  Using the field~$\field = \Z / 2\Z$, for~$b\neq \top$, $\BD_0[a,b]$ counts the number
    of zero-cycles in the simplicial complex~$\Ffunc(a)$ that are~zero-boundaries in the larger
    simplicial complex~$\Ffunc(b)$.  $\BD_0[a,\top]$ is a count of the number
    of vertices in~$\Ffunc(a)$.
    } \label{fig:ex-v-int2zb}
\end{figure}

\paragraph{Combinatorial PH Transform}
Now that we have a filtration functor, we can apply the persistent homology
pipeline of~\cite{PHFunc} and described in \secref{pipeline}.
Our next step is to apply the birth-death functor $\BD \colon \Fil \to \Mon$.
To keep the notation as simple as possible, we use the field $\field=\Z / 2\Z$.
We explain this functor by looking at how it behaves on the objects~$\Ffunc_{(+++)}$ and
$\Ffunc_{(0++)}$, and the arrow $\Ffunc_{(+++)} \to \Ffunc_{(0++)}$.
Recall from above that~$\Pcat_{(+++)}$ is the totally ordered lattice~$[v_1] <
[v_2] < [v_3] < [\top]$.
For simplicity of notation, we use~$[i,j]$ to represent~$\big[[v_i],[v_j]\big]$.
Thus,~$\Int{P_{(+++)}}$ is the nine element poset shown on the left hand side of
\subfigref{ex-v-int2zb}{interval}.
Similarly,~$\Int{P_{(0++)}}$ is the poset shown on the right hand side
of \subfigref{ex-v-int2zb}{interval}.
The bounded lattice function~$\alpha \colon \Pcat_{(+++)} \to \Pcat_{(0++)}$ (the
leftmost map in \figref{ex-v-lattice}) induces a bounded lattice function
$\intfcn{\alpha} \colon \Int{P_{(+++)}} \to \Int{P_{(0++)}}$.  In particular,
\begin{align*}
    [1,1]_{(0++)}
        & = \intfcn{\alpha}([1,1]_{(+++)})
        =\intfcn{\alpha}([1,2]_{(+++)})
        =\intfcn{\alpha}([2,2]_{(+++)}),\\
    [1,3]_{(0++)}
        & = \intfcn{\alpha}([1,3]_{(+++)})
        = \intfcn{\alpha}([2,3]_{(+++)}),\\
    [1,\top]_{(0++)}
        & = \intfcn{\alpha}([1,\top]_{(+++)})
        = \intfcn{\alpha}([2,\top]_{(+++)}),\\
    [3,3]_{(0++)}
        & = \intfcn{\alpha}([3,3]_{(+++)}),\\
    [3,\top]_{(0++)}
        & = \intfcn{\alpha}([3,\top]_{(+++)}), \text{ and }\\
    [\top,\top]_{(0++)}
        & = \intfcn{\alpha}([\top,\top]_{(+++)}).
\end{align*}
These three equivalences are illustrated with pink shading in \subfigref{ex-v-int2zb}{interval}.  If $i \neq 1$, then~$K$ has no
$i$-cycles and hence $\BD_i$ is trivial everywhere.  The values of $\BD_0$ on
the objects in~$\Int{P_{(+++)}}$ and $\Int{P_{0++}}$ are illustrated in
\subfigref{ex-v-int2zb}{zb}, with the
induced monotone-preserving morphism.

Finally, we apply the M\"obius Inversion Functor $\M \colon \Mon \to \Fnc$.
Again, we illustrate this functor by describing
how it behaves on $\BD_0(\Int{P_{(+++)}} \to \Int{P_{(0++)}})$.
First, we define~$\M_{(+++)} := \M (\BD_0(\Int{P_{(+++)}})$ to be the unique
function $\sigma \colon \Int{P_{(+++)}} \to \Z$ that satisfies \eqnref{MI},
where~$f=\BD_0(\Int{P_{(+++)}}$.  In particular, for $x \in \Int{P_{(+++)}}$
\[
    \M_{(+++)}(x)=
    \begin{cases}
        1 & x \in \big\{ [2,2], [3,3], [1,\top] \big\} \\
        0 & \text{otherwise}.
    \end{cases}
\]
In other words, the augmented persistence diagram (defined
in~\cite[Def.~1]{fasy2020faithful})
has three points:~$(2,2), (3,3)$, and
$(1,\infty)$.
This corresponds to the persistence diagram with one off-diagonal point, namely
$(1,\infty)$, which corresponds to the connected component born at the height
of~$v_1$.
Next, we consider the cell $S_{(0++)}$.  Let $x\in \Int{P_{(0++)}}$. Then,
$$
    \M_{(0++)}(x)=
    \begin{cases}
        1 & x \in \left\{ [1,1], [3,3], [1,\top] \right\} \\
        0 & \text{otherwise}.
    \end{cases}
$$
Moreover, the map $\M_{(+++)} \to \M_{(0++)}$ corresponds to mapping the points
in the following way: $[2,2]_{(+++)} \mapsto [1,1]_{(0++)}=[2,2]_{(0++)}$, $[3,3]_{(+++)} \mapsto[3,3]_{(0++)}$, and
$[1,T]_{(+++)} \mapsto [1,T]_{(0++)}$.  We note that while $[2,2]$ and $[1,T]$
have distinct birth heights in $\M_{(+++)}$, their images $[1,1]$ and $[1,T]$ have
the same birth height in $\M_{(0++)}$.

% BTF comment: removing this example
% \subsection{W}
% \todo{W example; can have a shorter exposition than the V}

\subsection{3D Example}

\begin{figure}
    \centering
    \begin{subfigure}[b]{0.4\textwidth}
        \centering
        \includegraphics[height=1.5in]{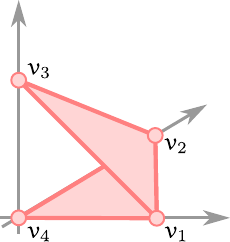}
        \caption{Complex in $\R^3$}
        \label{fig:ex-3d-gsc}
    \end{subfigure}
    \hspace{3em}
    \begin{subfigure}[b]{0.4\textwidth}
        \centering
        \includegraphics[height=1.5in]{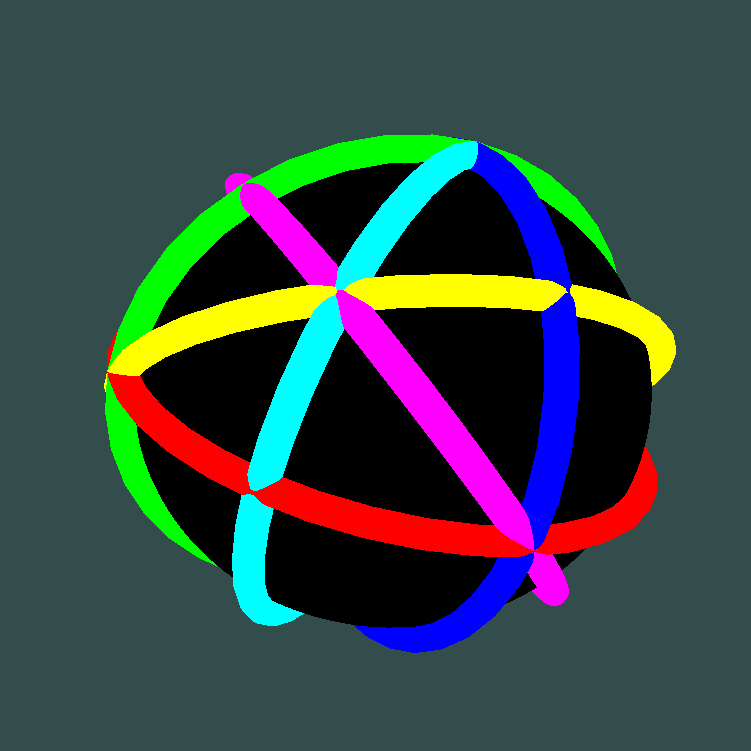}
        \caption{Cellulation of $\S^2$}
        \label{fig:ex-3d-strat}
    \end{subfigure}

    \caption{A geometric complex in $\R^3$ with
    the induced cellulation of $\S^2$. The cellulation depends only on the
    vetices of the complex.  Since the vertices are in general position, each
    great circle on~$\S^2$ is distinct.}
    \label{fig:ex-3d}
\end{figure}

We now consider a geometric complex embedded in $\R^3$.
Specifically, let $K$ be the following abstract simplicial complex known as a
book with two pages (see \subfigref{ex-3d}{gsc}):
$$
    K= \big\{ \{v_1\}, \{v_2\}, \{v_3\} , \{v_1,v_2\}, \{ v_2,v_3\},  \{ v_1,v_3\}, \{v_1, v_2, v_3\}
         \{v_1,v_4\}, \{v_2,v_4\}, \{v_1, v_2, v_4\} \big\}.
$$
This
complex is piecewise linearly embedded into $\R^3$ through the function $\phi
\colon |K| \to \R^{2}$ with~$\phi(v_1)=(1,0,0)$, $\phi(v_2)= (0,1,0)$, $\phi(v_3)=(0,0,1)$
and~$\phi(v_0)=(0,0,0)$. Note that the affine space spanned by
the points $\{v_1,v_2,v_3,v_4\}$ is $\R^3$, indicating that the vertices of~$K$ are in
general~position.
Rather than go through the whole pipeline again, we investigate the cell complex
and the poset of equivalence classes of vertices associated to each cell.

Let $(\S^2,\Sstrat)$ be the cell
complex defined in \propref{cellulation}.  This cellulation has $24$ two-cells,
one for each permutation of the vertices. (Note that if we were to add one more
vertex to $K$ that there would then be more permutations of the vertices than
cells of the cellulation due to geometric constraints). Each of the $36$ one-cells correspond to where
two vertices are at the same height and the remaining two vertices are at
distinct heights.  There are two types of zero-cells: (type-1) eight have degree six and
correspond to exactly three vertices at the same height, and (type-2) six have
degree four and correspond to two vertices at one height and the other two
vertices at another height.  See \subfigref{ex-3d}{strat}.

For each vector~$X \in \{-,0,+\}^6$, we use $\cell_{(X)}$ to denote the cell of
$\Sstrat$ labeled $X$ according to which side of $S_{1,2}$, $S_{1,3}$,
$S_{2,3}$, $S_{1,4}$, $S_{2,4}$, and $S_{3,4}$ the cell falls (if such a cell
exists). For example, $\cell_{(+---++)}=S_{1,2}^+ \cap S_{1,3}^- \cap S_{2,3}^- \cap
S_{1,4}^- \cap S_{2,4}^+ \cap S_{3,4}^+$.  The faces of $\cell_{(+---++)}$ are:
\begin{equation*}
    \{ \cell_{(+---++)}, \cell_{(0---++)}, \cell_{(+0--++)}, \cell_{(+---+0)},
        \cell_{(+0-0+0)}, \cell_{(0---+0)}, \cell_{(000-++)}
    \}.
\end{equation*}
In fact, every two-cell in this
cellulation has seven faces (itself, three one-cells, and three zero-cells).  In
other words, every two-cell is a topological triangle.  Some labelings of
zero-cells do not exist.  For example, no cell is labeled ${(00--++)}$.  By
contradiction, suppose there exists $\mu \in C_{(00--++)}$.  Then, since the
first coordinate is zero, $\mu \in
S_{1,2}$, which implies that $[v_1]_{\mu}=[v_2]_{\mu}$.  In addition, since the second
coordinate is zero, $\mu \in S_{2,3}$, which implies that $[v_2]_{\mu}=[v_3]_{\mu}$.
Thus, $[v_1]_{\mu}=[v_2]_{\mu}$, which means that $\mu \in S_{2,3}$, a contradiction.

To see what happens at a type-$1$ zero-cell,
consider $\Pcat_{(000-++)}$.
In particular, the equivalence classes for the vertices induced by a direction
vector in $C_{(000-++)}$ is:
\begin{align*}
    &  [v_1]_{(000-++)}  = [v_2]_{(000-++)} =  [v_3]_{(000-++)}\\
    &  [v_4]_{(000-++)}\\
    &  [\top].
\end{align*}
In general, all type-$1$ zero-cells have three vertices sharing an
equivalence class.
To see what happens at a type-$2$ zero-cell,
consider $\Pcat_{(0---+0)}$.
Here, we have the following equivalence~classess:
\begin{align*}
    &  [v_1]_{(0---+0)}  = [v_2]_{(0---+0)}\\
    &  [v_3]_{(0---+0)}  = [v_4]_{(0---+0)}\\
    &  [\top].
\end{align*}
In general, all type-$2$ zero-cells partition the vertices into two equivalence
classes, each of size two (and one additional equivalence class for $\top$).

\begin{figure}
    \centering
    \centering
    \includegraphics[height=2in]{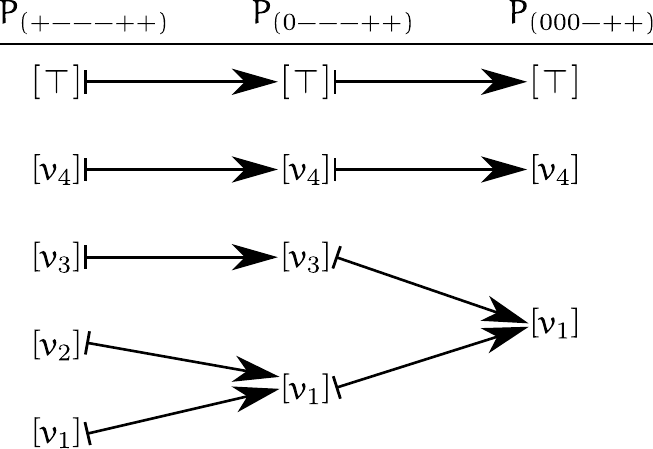}
    \caption{Face relations induce surjective poset maps.  These maps, in turn,
    induce arrows in~$\Fil$.
    } \label{fig:ex-3d-lattice}
\end{figure}
For geometric complexes in $\R^3$ or higher, we can look at compositions of
proper face relations.
Consider the following two face relations: $\cell_{(000-++)} < \cell_{(0---++)} <
\cell_{(+---++)}$.  These face relations induced posets of equivalence classes of the vertices, as
well as the maps connecting these posets, are given in
\figref{ex-3d-lattice}.
These two maps corresponds to the equivalences classes~$[v_1]_{(+---++)}\neq
[v_2]_{(+---++)}$ being mapped to the same equivalence class~$[v_1]_{(0---++)}$
in $P_{(0---++)}$, followed by the equivalences classes~$[v_1]_{(0---++)}\neq
[v_3]_{(0---++)}$
equivalence classes in~$P_{(+---++)}$
being mapped to the same equivalence class~$[v_1]_{(000-++)}$ in~$P_{(000-++)}$.
Observe that the
map~ $\Pcat_{(+---++)} \to \Pcat_{(000-++)}$
is exactly this composition resulting in collapsing three
equivalence classes in $P_{(+---++)}$ (namely,~$[v_1]$,~$[v_2]$, and~$[v_3]$)
to one equivalence~class.

%%%%%%%%%%%%%%%%%%%%%%%%%%%%%%%%%%%%%%%%%%%%%%

%%%%%%%%%%%%%%%%%%%%%%%%%%%%%%%%%%%%%%%%%%%%%%
\section{Discussion}\label{sec:discussion}

By taking a combinatorial approach to persistent homology, we are able to express
the PH transform of a geometric complex as a functor from the cellulation of a sphere to the category
of combinatorial persistence diagrams.
This point of view may be applied to any application that involves a parameterized family
of persistence diagrams.
For example, consider time-varying data \cite{frechet_time_varying, Kim_Memoli_2021, vineyard, time_varying_Morse}.
Assuming the data is finite, there is a cellulation of the real line (time parameter)
such that any two times on the same cell have the same combinatorial persistence diagram.
The resulting family of combinatorial persistence diagrams connected by charge-preserving morohisms
provides a categorical framework
in which to talk about all the persistence diagrams at once in a cohesive manner.

We hope this paper will start conversations on parameterized families of combinatorial persistence diagrams
in general and not just for the PH transform.
There are many questions to ask.
For example, consider a cell complex $(X, \Sstrat)$ on any space $X$ and a
functor~$F : \Sstrat \to \Fnc$.
What are its global co-sections?
Does the edit distance of McClearly and Patel lead to a stability result for such functors?
These are just a few questions worth asking.

%%%%%%%%%%%%%%%%%%%%%%%%%%%%%%%%%%%%%%%%%%%%%%

\paragraph{Acknowledgements}

The first author acknowledges Vin de Silva for a late night conversation
about counting cells and
David L. Millman for generating
\subfigref{ex-3d}{strat}.

\bibliographystyle{plain}
\bibliography{ref}{}

%%%%%%%%%%%%%%%%%%%%%%%%%%%%%%%%%%%%%%%%%%%%%%
\appendix
\FloatBarrier
\section{Notation and Definitions}\label{sec:background}

\paragraph{Posets}
Given a finite, totally ordered poset $P$, denote by $\bot \in P$ its smallest
element and by~$\top \in \Pcat$ the largest element.
A function $\alpha : \Pcat \to \Qcat$ between two finite, totally ordered posets
is a \emph{bounded monotone function} if $f(\top) = \top$, $f(\bot) = \bot$, and
for all $a \leq b$,~$\alpha(a) \leq \alpha(b)$.

For a finite, totally ordered poset $\Pcat$,
let $\Int{P}:= \big\{ [a,b] \subseteq \Pcat :
a \leq b \big \}$ be its set of non-empty intervals
with the following partial relation: $[a, b] \preceq [c, d]$ if $a \leq c$ and $b \leq d$.
Its bottom element is $[\bot, \bot]$ and its top element is~$[\top, \top]$.
A bounded monotone function $\alpha : \Pcat \to \Qcat$ between two finite, totally ordered posets induces a bounded monotone function $\bar \alpha : \bar P \to \bar Q$ where
$\bar \alpha \big( [a,b] \big) := \big[ \alpha(a), \alpha(b) \big]$; see \cite{PHFunc}
for a proof.

\paragraph{Categories and Functors}
In this paper,
we assume that the reader is familiar with categories and functors.
We defer to \cite{Riehl} for an introduction to category theory.
In Table~\ref{table:notation}, we provide notations used in this paper.

\begin{table}[htb]
    \caption{Notations for categories and functors. We assume
    $\Ccat$ is a category and $a,b,c \in \ob{\Ccat}$.}
    \label{table:notation}
    \centering
    \begin{tabular}{ | r | p{4in} | }
        \hline
        $\ob{\Ccat}$ & objects in $\Ccat$ \\ \hline
        $\Hom_\Ccat(a,b)$ & the set of morphisms or arrows between $a$ and $b$ in
        $\Ccat$ \\ \hline
        $\circ$ & composition of morphisms \\ \hline
        $1_a$ & the identity morphism in $\Hom_{\Ccat}(a,a)$ \\
        \hline
    \end{tabular}
\end{table}

\end{document}